\documentstyle{amsart}
\input{bull-art}
\newtheorem{thm}{Theorem}[section]

\newtheorem{cor}[thm]{Corollary}
\newtheorem{prop}[thm]{Proposition}
\newtheorem{conj}[thm]{Conjecture}
\newtheorem{flt}{Fermat's Last Theorem}

\newtheorem{ts}{Taniyama-Shimura Conjecture}

\newtheorem{ssts}{Semistable \hfil Taniyama-Shimura \hfil 
Conjecture}

\theoremstyle{definition}
\newtheorem{defn}{Definition}
\newtheorem{ex}{Example}
\newtheorem{exs}{\bf Examples}

\theoremstyle{remark}
\newtheorem{rem}{Remark}
\newtheorem{rems}{Remarks}
\newtheorem{notation}{Notation}

\newenvironment{mylist}{\begin{list}{$\bullet$}{
     \setlength{\leftmargin}{2.6 em}
     \setlength{\labelsep}{.4 em}
     \setlength{\labelwidth}{2.2 em}
     \setlength{\rightmargin}{0 em}}}{\end{list}}
\newenvironment{fermatquote}{\begin{list}{}{
     \setlength{\leftmargin}{1.5 em}
     \setlength{\rightmargin}{1.5 em}
     \setlength{\topsep}{.07 in}
     \setlength{\listparindent}{\parindent}}}{\end{list}}

\def\Q{{\bold Q}}
\def\Z{{\bold Z}}
\def\C{{\bold C}}
\def\O{{\cal O}}
\def\Of{{\cal O}_f}
\def\Ofl{{\cal O}_{f,\lambda}}
\def\D{{\cal D}}

\def\J{{\cal J}}
\def\L{{\cal L}}
\def\cc{{\cal C}}
\def\m{{\frak m}}
\def\q{{\cal Q}}
\def\a{{\frak a}}
\def\p{{\frak p}}
\def\P{{\cal P}}
\def\psubr{{\frak p}_R}
\def\psubt{{\frak p}_\T}
\def\l{{\lambda}}

\def\n{{\eta}}
\def\j{{\varphi}}
\def\e{{\varepsilon}}
\def\eps{{\epsilon}}
\def\wp{{{\bar \e}_p}}
\def\y{{\psi}}
\def\kd{{\theta_d}}

\def\H{{\frak H}}
\def\Mhat{{\hat M}}
\def\T{{\bold T}}
\def\Vn{V_n}
\def\Zp{{\bold Z}_p}
\def\Qp{{\bold Q}_p}
\def\Fp{{\bold F}_p}
\def\F{{\bold F}}
\def\Fq{{\bold F}_q}

\def\mn{{\boldsymbol \mu}_m}
\def\GQ{G_\Q}
\def\Gal{\mathrm{Gal}}
\def\Hom{\mathrm{Hom}}

\def\det{\mathrm{det}}
\def\im{\mathrm{Im}}
\def\trace{\mathrm{trace}}
\def\ker{\mathrm{ker}}
\def\res{\mathrm{res}}

\def\GL2{\mathrm{GL}_2}
\def\SL2{\mathrm{SL}_2}
\def\PGL2{\mathrm{PGL}_2}
\def\Sym2{\mathrm{Sym}^2}
\def\Frobq{\mathrm{Frob}_q}
\def\Frobmq{\mathrm{Frob}_\q}

\def\isom{\; {\buildrel \sim \over \rightarrow} \;}

\def\ABCD{\bigl(\begin{smallmatrix} a&b\\c&d 
\end{smallmatrix}\bigr)}
\def\fexp{\sum_{n=1}^{\infty}a_ne^{2 \pi inz}}
\def\sls{{S_\L^*}}
\def\Jqs{{J_q^*}}
\def\FLT{Fermat's Last Theorem }

\def\r{{\rho}}
\def\rb{{\bar\rho}}
\def\rhomod#1{{\rb_{E,#1}}}
\def\rhopmod#1{{\rb_{E',#1}}}
\def\rhoEp{{\r_{E,p}}}
\def\rhopEp{{\r_{E',p}}}
\def\rfl{{\rho_{f,\l}}}

\thicklines

\def\downdiag#1#2#3#4#5{
   \begin{center}
   \begin{picture}(85,50)
   \put(5,40){\makebox(0,0){$\GQ$}}
   \put(80,40){\makebox(0,0){$\GL2(#1)$}}
   \put(80,5){\makebox(0,0){$\GL2(#2)$}}
   \put(18,40){\vector(1,0){38}}
   \put(80,30){\vector(0,-1){15}}
   \put(18,33){\vector(2,-1){41}}
   \put(37,44){\makebox(0,0)[b]{\scriptsize $[#3]$}}
   \put(36,18){\makebox(0,0)[r]{\scriptsize $[#4]$}}
   \put(83,24){\makebox(0,0)[l]{\scriptsize $#5$}}
   \end{picture}
   \end{center}}

\def\updiag#1#2#3#4#5{
   \begin{center}
   \begin{picture}(85,55)
   \put(5,15){\makebox(0,0){$\GQ$}}
   \put(80,15){\makebox(0,0){$\GL2(#2)$}}
   \put(80,50){\makebox(0,0){$\GL2(#1)$}}
   \put(18,15){\vector(1,0){40}}
   \put(80,40){\vector(0,-1){15}}
   \put(18,26){\vector(2,1){37}}
   \put(38,12){\makebox(0,0)[t]{\scriptsize $[#4]$}}
   \put(32,38){\makebox(0,0)[br]{\scriptsize $[#3]$}}
   \put(83,34){\makebox(0,0)[l]{\scriptsize $#5$}}
   \end{picture}
   \end{center}}

\title{A report on Wiles' Cambridge lectures}
\author[K. Rubin]{K. Rubin}
\address{Department of Mathematics, Ohio State University, 
Columbus, Ohio 43210}
\email{rubin\char`\@math.ohio-state.edu}
\author[A. Silverberg]{A. Silverberg}
\address{Department of Mathematics, Ohio State University, 
Columbus, Ohio 43210}
\email{silver\char`\@math.ohio-state.edu}
\date{November 29, 1993}
\thanks{The authors thank the National Science Foundation 
 for financial support}
\begin{document}
\def\currentvolume{31}
\def\currentissue{1}
\def\currentyear{1994}
\def\currentmonth{July}
\def\copyrightyear{1994}
\def\currentpages{15-38}

\maketitle
\begin{abstract}
In lectures at the Newton Institute in June of 1993, 
Andrew Wiles announced
a proof of a large part of the Taniyama-Shimura Conjecture 
and, as a consequence,
Fermat's Last Theorem.  This report for nonexperts 
discusses the mathematics 
involved in Wiles' lectures, including the necessary 
background and the
mathematical history.
\end{abstract}

\section*{Introduction}  

On June 23, 1993, Andrew Wiles wrote on a blackboard, 
before an 
audience at the Newton Institute in Cambridge, England, 
that if $p$
is a prime number, $u$, $v$, and $w$ are rational numbers, 
and 
$u^p + v^p + w^p = 0$, then $uvw = 0$.  In other words, he 
announced that 
he could prove Fermat's Last Theorem.  His announcement 
came at the end 
of his series of three talks entitled ``Modular forms, 
elliptic curves, and Galois 
representations'' at the week-long workshop on 
``$p$-adic Galois representations, Iwasawa theory, and 
the Tamagawa numbers of motives''.

In the margin of his copy of the works of Diophantus, next 
to a problem on 
Pythagorean triples, Pierre de Fermat (1601--1665) wrote: 
\begin{fermatquote}
\item \hskip\listparindent 
{\em Cubum autem in duos cubos, aut quadratoquadratum in 
duos quadratoquadratos, et 
generaliter nullam in infinitum ultra quadratum potestatem 
in duos ejusdem 
nominis fas est dividere\,\RM: cujus rei demonstrationem 
mirabilem sane detexi. Hanc 
marginis exiguitas non caperet.}  

(It is impossible to separate a cube into two cubes, 
or a fourth power into two fourth powers, or in general, 
any power higher than the 
second into two like powers.  I have discovered a truly 
marvelous proof of this, 
which this margin is too narrow to contain.)
\end{fermatquote}
We restate Fermat's conjecture as follows.

\begin{flt}
If $n > 2$, then $a^n + b^n = c^n$ has no solutions in 
nonzero 
integers $a$, $b$, and $c$.
\end{flt}

A proof by Fermat has never been found,
 and the problem has remained open, inspiring
many generations of mathematicians. Much of modern number 
theory has been
built on attempts to prove Fermat's Last Theorem.  For 
details on the history 
of Fermat's Last Theorem (last because it is the last of 
Fermat's questions to be 
answered) see \cite{Dickson}, \cite{Edwards}, and 
\cite{Ribenboim}.

What Andrew Wiles announced in Cambridge was that he could 
prove ``many'' 
elliptic curves  are modular, sufficiently many to imply 
Fermat's Last Theorem.  In 
this paper we will explain Wiles' work on elliptic curves 
and its connection 
with Fermat's Last Theorem.  In \S\ref{ellipticcurves} we 
introduce elliptic 
curves and modularity, and give the connection 
between  Fermat's Last Theorem and the Taniyama-Shimura 
Conjecture on the modularity 
of elliptic curves.  
In \S\ref{overview} we describe how Wiles reduces the 
proof of the 
Taniyama-Shimura Conjecture to what we call the Modular 
Lifting 
Conjecture (which can be viewed as a weak form of the
Taniyama-Shimura Conjecture), by using a theorem of 
Langlands and Tunnell.
In \S\ref{representations}  
and \S\ref{universal} we show how the Semistable Modular 
Lifting Conjecture 
is related to a conjecture of Mazur on deformations 
of Galois representations (Conjecture \ref{M}), and in 
\S\ref{proof} we
describe Wiles' 
method of attack on this conjecture.  In order to make 
this survey as accessible 
as possible to nonspecialists, the more technical details 
are postponed as 
long as possible, some of them to the appendices.

Much of this report is based on Wiles' lectures in 
Cambridge. The authors
apologize for any errors we may have introduced. We also 
apologize to those
whose mathematical contributions we, due to our incomplete 
understanding, do not
properly acknowledge.

The ideas Wiles introduced in his Cambridge lectures will 
have an important
influence on research in number theory. Because of the 
great interest 
in this subject and the lack of a publicly available 
manuscript, we hope this 
report will be useful to the mathematics community. In 
early December 1993,
shortly before this paper went to press, Wiles announced 
that ``the final
calculation of a precise upper bound for the Selmer group 
in the semistable
case'' (see \S5.3 and \S5.4 below) ``is not yet complete 
as it stands,'' but
that he believes he will be able to finish it in the near 
future using the
ideas explained in his Cambridge lectures. While Wiles' 
proof of Theorem \ref{BW}
below and Fermat's Last Theorem depends on the calculation 
he referred to in his 
December announcement, Theorem \ref{LW} and Corollary 
\ref{cm} do not. Wiles' work 
provides for the first time infinitely many modular 
elliptic curves over the 
rational numbers which are not isomorphic over the complex 
numbers (see
\S\ref{fr} for an explicit infinite family).

\begin{notation}
The integers, rational numbers, complex numbers, and 
$p$-adic integers
will be denoted $\Z$, $\Q$, $\C$, and $\Z_p$, respectively. 
If $F$ is a field, then ${\bar F}$ denotes an
algebraic closure of $F$.
\end{notation}


\section{Connection between \FLT and elliptic curves}
\label{ellipticcurves}

\subsection{\FLT follows from modularity of elliptic curves}
\label{mec}

Suppose \FLT were false.  Then there would exist nonzero 
integers $a$, $b$, $c$, and  $n > 2$ such that $a^n + b^n 
= c^n$.  It is
easy to see that no generality is lost by assuming that 
$n$ is a prime
greater than three (or greater than four million, by 
\cite{Buhler}; 
see \cite{Hardy-Wright} for $n = 3$ and $4$) and that
$a$ and $b$ are relatively prime.  Write down the cubic 
curve:
\begin{equation}
\label{ell}
	y^2 = x(x + a^n)(x - b^n).
\end{equation}

In \S\ref{E} we will see that such curves are elliptic 
curves, and in
\S\ref{m1} we will explain what it means for an elliptic 
curve to be
modular. Kenneth Ribet \cite{Ribet} proved that if $n$ is 
a prime greater 
than three, $a$, $b$, and $c$ are nonzero integers,
and $a^n + b^n = c^n$, then the elliptic curve 
(\ref{ell}) 
is not modular.  But the results announced by Wiles imply 
the following.

\begin{thm}[Wiles]
\label{W}
If $A$ and $B$ are distinct, nonzero, relatively prime 
integers, and 
$AB(A-B)$ is divisible by $16$, then the elliptic curve 
$$
y^2 = x(x + A)(x + B)
$$
is modular.
\end{thm}

Taking $A = a^n$ and $B = -b^n$ with $a$, $b$, $c$, and 
$n$ coming from
our hypothetical solution to a Fermat equation as above, 
we see that
the conditions of Theorem \ref{W} are satisfied since $n 
\geq 5$ and
one of $a$, $b$, and $c$ is even.  Thus Theorem \ref{W} 
and Ribet's result
together imply Fermat's Last Theorem!

\subsection{History}
The story of the connection between \FLT and elliptic 
curves begins in
1955, when Yutaka Taniyama (1927--1958) posed problems 
which may be
viewed as a weaker version of
the following conjecture (see \cite{Shimura-obit}).

\begin{ts}
Every elliptic curve over $\Q$ is modular.
\end{ts}

The conjecture in the present form was made by Goro 
Shimura around 1962--64 
and has become better understood due to work of Shimura 
[33--37]
and of Andr\'e Weil \cite{Weil} (see also \cite{Eichler}). 
The 
Taniyama-Shimura Conjecture is one of the major 
conjectures in number theory.

Beginning in the late 1960s 
[15--18], Yves Hellegouarch 
connected Fermat equations $a^n + b^n = c^n$ with elliptic 
curves of the form
(1) and used results about \FLT
to prove results about elliptic curves.  The landscape 
changed abruptly
in 1985 when Gerhard Frey stated in a lecture at 
Oberwolfach that
elliptic curves arising from counterexamples to \FLT could 
not be modular
\cite{Frey}.  Shortly thereafter Ribet \cite{Ribet} proved 
this, following ideas of 
Jean-Pierre Serre \cite{Serre} (see \cite{Oesterle} for a 
survey). In other words,
``Taniyama-Shimura Conjecture $\Rightarrow$ Fermat's Last 
Theorem''.  

Thus, the stage was set. A proof of the Taniyama-Shimura 
Conjecture
(or enough of it to know that elliptic curves coming 
from Fermat equations are modular)
would be a proof of Fermat's Last Theorem.

\subsection{Elliptic curves}  
\label{E}

\begin{defn}
An {\em elliptic curve} over $\Q$ is a nonsingular curve 
defined by an
equation of the form
\begin{equation}
\label{wf}
		y^2 + a_1xy + a_3y  = x^3 + a_2x^2 + a_4x + a_6
\end{equation}
where the coefficients $a_i$ are integers.  The solution
$(\infty,\infty)$ will be viewed as a point on the 
elliptic curve.
\end{defn}

\begin{rems}
(i) A {\em singular point} on a curve $f(x,y) = 0$ is a 
point where 
both partial derivatives vanish.  A curve is {\em 
nonsingular} if it has no
singular points.

(ii) Two elliptic curves over $\Q$ are {\em
isomorphic} if one can be obtained from the other by 
changing
coordinates $x = A^2x' + B$, $y = A^3y' + Cx' + D$, with 
$A$, $B$, $C$,
$D \in \Q$ and dividing through by $A^6$.

(iii) Every elliptic curve over $\Q$ is isomorphic to one 
of the
form
$$
		y^2  = x^3 + a_2x^2 + a_4x + a_6
$$
with integers $a_i$.  A curve of this form is nonsingular 
if and 
only if the cubic on the right side has no repeated roots.
\end{rems}

\begin{ex}
The equation $y^2 = x(x + 3^2)(x - 4^2)$ defines an 
elliptic curve
over $\Q$.
\end{ex}

\subsection{Modularity}  
\label{m1}
Let $\H$ denote the complex upper half plane 
$\{z \in \C : \im(z) > 0\}$ where $\im(z)$ is the 
imaginary part of $z$.  
If $N$ is a positive integer, define a group of matrices
$$
\Gamma_0(N) = \bigl\{\ABCD \in 
\operatorname{SL}_2(\Z) : {\text{$c$ is divisible by $N$}} 
\bigr\}.
$$
The group $\Gamma_0(N)$ acts on $\H$ by linear fractional 
transformations
$
\ABCD(z)= {{az+b} \over {cz+d}}.
$
The quotient space $\H/\Gamma_0(N)$ is a (noncompact) 
Riemann
surface.  It can be completed to a compact Riemann 
surface, denoted
$X_0(N)$, by adjoining a finite set of points called 
cusps. The cusps are 
the finitely many
equivalence classes of $\Q \cup \{i\infty\}$ under the 
action of
$\Gamma_0(N)$ (see Chapter 1 of \cite{Shimura-red-book}).  
The complex points of an elliptic curve can also be viewed 
as a 
compact Riemann surface.

\begin{defn}
An elliptic curve $E$ is {\em modular} if, for some 
integer $N$, there
is a holomorphic map from $X_0(N)$ onto $E$.
\end{defn}

\begin{ex}
It can be shown that there is a (holomorphic) isomorphism 
from $X_0(15)$ 
onto the elliptic curve 
$y^2 = x(x + 3^2)(x - 4^2)$.
\end{ex}

\begin{rem}
There are many equivalent definitions of modularity (see 
{\S}II.4.D of 
\cite{Oesterle} and appendix of \cite{Mazur}). In some 
cases the equivalence
is a deep result. For Wiles' proof of Fermat's Last 
Theorem it suffices to
use only the definition given in \S \ref{mdeftwo} below.
\end{rem}

\subsection{Semistability}
\label{ss}

\begin{defn}
An elliptic curve over $\Q$ is {\em semistable at the 
prime} $q$ if it is 
isomorphic to an elliptic curve over $\Q$
which modulo $q$ either is nonsingular or has a singular 
point with
two distinct tangent directions.  An elliptic curve over 
$\Q$ is called 
{\em semistable} if it is semistable at every prime.
\end{defn}

\begin{ex}
The elliptic curve $y^2 = x(x + 3^2)(x - 4^2)$ is 
semistable because it
is isomorphic to $y^2+xy+y = x^3+x^2-10x-10$, but the 
elliptic 
curve $y^2 = x(x + 4^2)(x - 3^2)$ is not semistable (it is 
not semistable at 2).
\end{ex}

Beginning in \S\ref{overview} we explain how Wiles shows 
that his main result 
on Galois representations (Theorem \ref{BW}) implies the 
following part of
the Taniyama-Shimura Conjecture.

\begin{ssts}
\label{SS}
Every \hfil semistable \hfil elliptic \hfil curve over 
$\Q$ is modular.
\end{ssts}

\begin{prop}
\label{stscflt}
The Semistable Taniyama-Shimura Conjecture implies Theorem 
\ref{W}.
\end{prop}

\begin{pf}
If  $A$  and  $B$  are distinct, nonzero, relatively prime 
integers, write 
$E_{A,B}$ for the elliptic curve defined by  $y^2 = x(x+
A)(x+B)$.  
Since $E_{A,B}$ and $E_{-A,-B}$ are 
isomorphic over the complex numbers (i.e., as Riemann 
surfaces), $E_{A,B}$ is 
modular if and only if $E_{-A,-B}$ is modular.  If further 
$AB(A-B)$ is divisible by  16, then either $E_{A,B}$ or 
$E_{-A,-B}$ is 
semistable (this is easy to check directly; see for 
example {\S}I.1
of \cite{Oesterle}).  The Semistable Taniyama-Shimura 
Conjecture now implies 
that both $E_{A,B}$ and $E_{-A,-B}$ are modular, and thus 
implies Theorem \ref{W}.
\end{pf}

\begin{rem}
In \S\ref{mec} we saw that Theorem \ref{W} and Ribet's 
Theorem together
imply Fermat's Last Theorem. Therefore, the Semistable 
Taniyama-Shimura 
Conjecture implies Fermat's Last Theorem.
\end{rem}

\subsection{Modular forms}  
\label{modforms}
In this paper we will work with a definition of 
modularity which uses modular forms. 

\begin{defn}
\label{mf}
If $N$ is a positive integer, a {\em modular form} $f$ of
weight $k$ for $\Gamma_0(N)$ is a holomorphic function $f 
: \H \to \C$ which
satisfies 
\begin{equation}
\label{transf}
f(\gamma(z)) = (cz+d)^kf(z)  {\text{\quad for every~}} 
\gamma = 
\bigl( \begin{smallmatrix} a & b \\ c & d 
\end{smallmatrix} \bigr) 
\in \Gamma_0(N){\text{~and~}} z \in \H,
\end{equation}
and is holomorphic at the cusps (see Chapter 2 of 
\cite{Shimura-red-book}).  
\end{defn}

A modular form $f$ satisfies $f(z) = f(z+1)$ (apply 
(\ref{transf})  to
$\bigl({\begin{smallmatrix} 1&1\\0&1 
\end{smallmatrix}}\bigr) \in \Gamma_0(N)$),
so it has a Fourier expansion
$f(z) = \sum_{n=0}^{\infty}a_ne^{2 \pi inz}$,  with 
complex numbers 
$a_n$ and with $n \geq 0$
because $f$ is holomorphic at the cusp $i\infty$.
We say $f$ is a {\em cusp form} if it vanishes at all the 
cusps; 
in particular for a cusp form the coefficient $a_0$ (the 
value at
$i\infty$) is zero.  Call a cusp form {\em normalized} if 
$a_1 = 1$.

 For fixed $N$ there are commuting linear operators 
(called {\em Hecke operators}) 
$T_m$, for integers $m \geq 1$, on the 
(finite-dimensional) vector space of 
cusp forms of weight two for $\Gamma_0(N)$
(see Chapter 3 of \cite{Shimura-red-book}).  If $f(z) = 
\fexp$, then
\begin{equation}
\label{tm}
T_mf(z) = \sum_{n=1}^\infty \bigl(\sum_{{(d,N)=1}\atop {d 
\mid (n,m)}}
da_{{nm}/{d^2}}\bigr) e^{2\pi inz}
\end{equation}
where $(a,b)$ denotes the greatest common divisor of $a$ 
and $b$ and $a\mid b$
means that $a$ divides $b$. 
The {\em Hecke algebra} $T(N)$ is the ring generated over 
$\Z$ by these operators.
	 
\begin{defn}
In this paper an {\em eigenform} will mean a normalized 
cusp form of weight two 
for some $\Gamma_0(N)$ which is an eigenfunction for all 
the Hecke operators.  
\end{defn}

By (\ref{tm}), if $f(z) = \fexp$ is an eigenform, 
then $T_mf = a_mf$ for all $m$.

\subsection{Modularity, revisited}
\label{mdeftwo}
Suppose  $E$  is an elliptic curve
over  $\Q$.  If $p$ is a prime, write $\Fp$ for the finite 
field with
$p$ elements, and let $E(\Fp)$ denote the $\Fp$-solutions 
of the equation
for $E$ (including the point at infinity).  
We now give a second definition of modularity for an 
elliptic curve.  

\begin{defn}
\label{m2}
An elliptic curve $E$ over $\Q$ is {\em modular} if there 
exists an
eigenform  $\fexp$ such that for all but finitely many 
primes $q$,
\begin{equation}
\label{mdlr}
a_q  =  q + 1 - \#(E(\Fq)).
\end{equation}
\end{defn}


\section{An overview}
\label{overview}

The flow chart shows how \FLT would follow if one knew
the Semistable Modular Lifting Conjecture (Conjecture 
\ref{ME}) for the
primes 3 and 5. In \S\ref{ellipticcurves} we discussed the 
upper arrow, i.e., 
the implication
``Semistable Taniyama-Shimura Conjecture $\Rightarrow$ 
Fermat's Last Theorem''. 
In this section we will discuss the other implications in 
the flow chart. 
The implication given by the lowest arrow is 
straightforward 
(Proposition \ref{irredat3}), while the middle one uses an 
ingenious
idea of Wiles (Proposition \ref{redat3}). 


\begin{figure}[hbt]
\begin{center}
\begin{picture}(305,250)

\put(120,210){\framebox(130,30){Fermat's Last Theorem}}
\put(185,185){\vector(0,1){20}}
\put(110,140){\framebox(150,40){\parbox{2in}{\begin{center} 
     Semistable Taniyama-Shimura \\ Conjecture 
\end{center}}}}
\put(185,115){\vector(0,1){20}}
\put(170,100){\line(1,1){15}}
\put(200,100){\line(-1,1){15}}
\put(15,70){\framebox(150,40){\parbox{2in}{\begin{center} 
     Semistable Taniyama-Shimura \\ for $\rhomod{3}$ 
irreducible \end{center}}}}
\put(205,70){\framebox(105,40){\parbox{2in}{\begin{center} 
     Semistable Modular \\ Lifting for $p = 5$ 
\end{center}}}}
\put(120,45){\vector(0,1){20}}
\put(105,30){\line(1,1){15}}
\put(135,30){\line(-1,1){15}}
\put(0,0){\framebox(100,40){\parbox{2in}{\begin{center} 
     Langlands-Tunnell \\ Theorem \end{center}}}}
\put(140,0){\framebox(105,40){\parbox{2in}{\begin{center} 
 Semistable Modular \\ Lifting for $p = 3$ \end{center}}}}
\end{picture}

\bigskip
Semistable Modular Lifting Conjecture $\Rightarrow$ \FLT.
\end{center}
\end{figure}


\begin{rem}
By the Modular Lifting Conjecture we will mean the 
Semistable Modular
Lifting Conjecture with the hypothesis of semistability 
removed. The
arguments of this section can also be used to show that 
the Modular
Lifting Conjecture for $p = 3$ and $5$, together with the 
Langlands-Tunnell Theorem, imply the full Taniyama-Shimura 
Conjecture.
\end{rem}

\subsection{Semistable Modular Lifting}
\label{mcec}
Let $\bar\Q$ denote the algebraic closure of $\Q$ in $\C$,
and let $\GQ$ be the Galois group $\Gal(\bar\Q/\Q)$.
If $p$ is a prime, write
$$
\wp : \GQ \to \Fp^\times
$$
for the character giving the action of $\GQ$ on the $p$-th 
roots of unity.  For 
the facts about elliptic curves stated below, see 
\cite{Silverman}.  
If $E$ is an elliptic curve over $\Q$ and $F$ is a 
subfield of the complex numbers, 
there is a natural commutative group law on the set of 
$F$-solutions of $E$, with the 
point at infinity as the identity element.  Denote this 
group $E(F)$.  
If $p$ is a prime, write  $E[p]$  for the subgroup of 
points in $E(\bar\Q)$ of order 
dividing  $p$. Then  $E[p] \cong \Fp^2$.  The action of  
$\GQ$  on  $E[p]$  
gives a continuous representation  
$$
\rhomod{p} : \GQ  \to  \operatorname{GL}_2(\Fp)
$$
(defined up to isomorphism) such that
\begin{equation}
\label{cycl}
\det(\rhomod{p}) = \wp
\end{equation}
and for all but finitely many primes $q$,
\begin{equation}
\label{es}
\trace(\rhomod{p}(\Frobq)) \equiv q + 1 - \#(E(\Fq)) 
\quad(\text{\rom{mod}}\ p).
\end{equation}
(See Appendix \ref{A1} for the definition of the Frobenius 
elements $\Frobq \in \GQ$ 
attached to each prime number $q$.)

If  $f(z) = \fexp$ is an
eigenform, let  $\Of$ denote the ring of integers of the 
number field
$\Q(a_2, a_3, \ldots)$. (Recall that our eigenforms are 
normalized so
that $a_1 = 1$.)

The following conjecture is in the spirit of a conjecture 
of Mazur (see
Conjectures \ref{M1} and \ref{M}).

\begin{conj}[Semistable \hfil Modular \hfil Lifting \hfil 
Conjecture]
\label{ME}
Suppose \hfil $p$ \hfil is \hfil an \hfil odd prime and 
$E$ is a semistable 
elliptic curve over $\Q$ satisfying

\rom{(a)} $\rhomod{p}$  is irreducible,

\rom{(b)} there are an eigenform $f(z) = \fexp$ and a 
prime ideal  $\lambda$  of $\Of$ such that 
$p \in \lambda$ and for all but finitely many primes $q$,
$$
a_q  \equiv  q + 1 - \#(E(\Fq))   \pmod{\lambda}.
$$
Then  $E$  is modular.
\end{conj}

The Semistable Modular Lifting Conjecture is a priori 
weaker than the 
Semistable Taniyama-Shimura Conjecture because 
of the extra hypotheses (a) and (b).  The more serious 
condition is (b); 
there is no known way to produce such a form in general.  
But when  
$p = 3$,  the existence of such a form follows from the 
theorem below of 
Tunnell \cite{Tunnell} and Langlands \cite{Langlands}.  
Wiles 
then gets around condition (a) by a clever argument 
(described below) 
which, when  $\rhomod{3}$  is not irreducible, allows him 
to use  $p = 5$  instead.

\subsection{Langlands-Tunnell \hfil Theorem}
In \hfil order \hfil to \hfil state \hfil the \hfil 
Langlands-Tunnell \hfil 
Theorem, we need weight-one modular
forms for a subgroup of $\Gamma_0(N)$.  Let
$$
\Gamma_1(N) = \bigl\{\ABCD \in 
\operatorname{SL}_2(\Z) : c \equiv 0 \!\pmod{N},\; 
a \equiv d \equiv 1 \!\pmod{N}\bigr\}.
$$
Replacing $\Gamma_0(N)$ by $\Gamma_1(N)$ in \S 
\ref{modforms}, one can
define the notion of cusp forms on $\Gamma_1(N)$. See 
Chapter 3 of 
\cite{Shimura-red-book} for the definitions of the Hecke 
operators on the
space of weight-one cusp forms for $\Gamma_1(N)$.

\begin{thm}[Langlands-Tunnell]
\label{LT}
Suppose  $\rho : \GQ  \to  \GL2(\C)$  is a continuous 
irreducible
representation whose image in  $\PGL2(\C)$  is 
a subgroup of  $S_4$ 
\RM(the symmetric group on four elements\,\RM),  $\tau$ is
complex conjugation, and $\det(\rho(\tau)) = -1$. Then 
there is a weight-one   
cusp form $\sum_{n=1}^{\infty}b_ne^{2 \pi inz}$ for some 
$\Gamma_1(N)$, which
is an eigenfunction for all the corresponding Hecke 
operators, such that for 
all but finitely many primes  $q$,
\begin{equation}
\label{LTform}
b_q  =  \trace(\rho(\Frobq)).
\end{equation}
\end{thm}

The theorem as stated by Langlands \cite{Langlands} 
and by Tunnell \cite{Tunnell} produces an automorphic 
representation rather than a cusp form.  Using the fact 
that 
$\det(\rho(\tau)) = -1$, standard techniques (see for 
example 
\cite{Gelbart}) show that this automorphic representation 
corresponds to a 
weight-one cusp form as in Theorem \ref{LT}.

\subsection{Semistable Modular Lifting $\Rightarrow$ 
Semistable Taniyama-Shimura}

\begin{prop}
\label{irredat3}
Suppose the Semistable Modular Lifting Conjecture is true 
for  $p = 3$, $E$  is a 
semistable elliptic curve, and  $\rhomod{3}$  is 
irreducible.  Then  $E$  is modular.
\end{prop}

\begin{pf}
It suffices to show that hypothesis (b) of the Semistable 
Modular Lifting Conjecture 
is satisfied with the given curve  $E$, for  $p = 3$.  
There is a  faithful representation
$$
\psi : \GL2(\F_3)  \hookrightarrow  \GL2(\Z[\sqrt{-2}])  
\subset  \GL2(\C)
$$
such that for every $g \in \GL2(\F_3)$,
\begin{equation}
\label{trace}
\trace(\psi(g)) \equiv \trace(g)\quad  (\text{\rom{mod}}(1+
\sqrt{-2}))
\end{equation}
and
\begin{equation}
\label{det}
\det(\psi(g)) \equiv \det(g)  \quad (\text{\rom{mod}}\ 3).
\end{equation}
Explicitly, $\psi$ can be defined on generators of 
$\GL2(\F_3)$ by 
$$\psi
\left(\left(\begin{array}{rr}-1&1\\-1&0\end{array}\right)%
\right) = 
\left(\begin{array}{rr}-1&1\\-1&0\end{array}\right)
\quad \text{and} \quad 
\psi\left(\left(\begin{array}{rr}1&-1\\1&1\end{array}%
\right)\right) = 
\left(\begin{array}{rr}\sqrt{-2}&1\\1&0\end{array}\right).$$
Let $\rho = \psi \circ \rhomod{3}$.
If $\tau$ is complex conjugation,
then it follows from (\ref{cycl}) and (\ref{det}) that 
$\det(\rho(\tau)) = -1$.
The image of  $\psi$  in  $\PGL2(\C)$  is a subgroup of 
$\PGL2(\F_3) \cong S_4$.  
Using that $\rhomod{3}$ is irreducible, one can show that 
$\r$ is irreducible.

Let $\p$ be a prime of $\bar\Q$ containing $1+\sqrt{-2}$.
Let $g(z) = \sum_{n=1}^{\infty}b_ne^{2 \pi inz}$ be a 
weight-one cusp form
for some $\Gamma_1(N)$ obtained by applying the 
Langlands-Tunnell Theorem 
(Theorem \ref{LT}) to $\rho$.  It follows from 
(\ref{cycl}) and (\ref{det})
that $N$ is divisible by $3$.
The function
$$
{\bold E}(z) = 1 + 6\sum_{n=1}^\infty \sum_{d \mid n} 
\chi(d)e^{2\pi inz}
\quad{\text{where}\ }
 \chi(d) = \left\{
\begin{array}{rl}
0 & \text{if $d \equiv 0 \pmod{3}$,} \\
1 & \text{if $d \equiv 1 \pmod{3}$,} \\
-1 & \text{if $d \equiv 2 \pmod{3}$}
\end{array}
\right.
$$
is a weight-one modular form for $\Gamma_1(3)$. The product 
$g(z){\bold E}(z) = \sum_{n=1}^{\infty}c_ne^{2 \pi inz}$
is a weight-two cusp form for $\Gamma_0(N)$ with $c_n 
\equiv b_n \pmod{\p}$ for
all $n$.  It is now possible to find an 
eigenform $f(z) = \fexp$ on $\Gamma_0(N)$
such that $a_n  \equiv  b_n \pmod{\p}$ for every  $n$ 
(see 6.10 and 6.11 of \cite{Deligne-Serre}).
By (\ref{es}), (\ref{LTform}), and (\ref{trace}),  
$f$  satisfies (b) of the Semistable Modular Lifting 
Conjecture with $p = 3$ and
with $\lambda = \p \cap \Of$.
\end{pf}

\begin{prop}[Wiles]
\label{redat3}
Suppose the Semistable Modular Lifting Conjecture is true 
for  $p = 3$  and  $5$,  
$E$  is a semistable elliptic curve 
over $\Q$, and  $\rhomod{3}$  is reducible.  Then  $E$  is 
modular.
\end{prop}

\begin{pf}
The elliptic curves over $\Q$ for which both  $\rhomod{3}$ 
 and  
$\rhomod{5}$  are reducible are all known to be modular 
(see Appendix \ref{app1}).  
Thus we can suppose  $\rhomod{5}$  is irreducible.  It 
suffices to produce an 
eigenform as in (b) of the Semistable Modular Lifting 
Conjecture,
but this time there is no analogue of the Langlands-Tunnell
Theorem to help.  Wiles uses the Hilbert Irreducibility 
Theorem, applied 
to a parameter space of elliptic curves, 
to produce another semistable elliptic curve  $E'$  over 
$\Q$ satisfying
\begin{mylist}
\item[(i)]  $\rhopmod{5}$  is isomorphic to  $\rhomod{5}$, 
and
\item[(ii)]  $\rhopmod{3}$  is irreducible.
\end{mylist}
(In fact there will be infinitely many such  $E'$; see 
Appendix \ref{app2}.)  
Now by Proposition \ref{irredat3},  $E'$  is modular.  Let  
$f(z) = \fexp$  be a corresponding eigenform.  Then for all 
but finitely many primes $q$,
\begin{equation*}
\begin{split}
a_q  = &\; q + 1 - \#(E'(\Fq))  \equiv  
\trace(\rhopmod{5}(\Frobq))  \cr
 \equiv& \; \trace(\rhomod{5}(\Frobq))  \equiv  q + 1 - 
\#(E(\Fq))  \quad
(\text{\rom{mod}}\ 5)
\end{split}
\end{equation*}
by (\ref{es}).  Thus the form  $f$  satisfies hypothesis 
(b) of  
the Semistable Modular Lifting Conjecture, and we conclude 
that  $E$  is modular.
\end{pf}

Taken together, Propositions \ref{irredat3} and 
\ref{redat3} show that 
the Semistable Modular Lifting Conjecture
for  $p = 3$  and  $5$  implies the Semistable 
Taniyama-Shimura Conjecture.


\section{Galois representations}
\label{representations}
The next step is to translate the Semistable Modular 
Lifting Conjecture into a
conjecture (Conjecture \ref{M1}) about the modularity of 
liftings of
Galois representations. Throughout this paper, if $A$ is a 
topological ring,
a representation $\rho : \GQ \to \operatorname{GL}
_2(A)$ will mean a continuous homomorphism
and $[\r]$ will denote the isomorphism class of $\r$.
If $p$ is a prime, let
$$\e_p : \GQ \to \Zp^\times$$  
be the character giving the action of $\GQ$ on $p$-power 
roots of unity.

\subsection{The $p$-adic representation attached to an 
elliptic curve}
Suppose  $E$  is an elliptic curve over  $\Q$ and $p$ is a 
prime number.  For 
every positive integer $n$, write  $E[p^n]$  for the 
subgroup  in  $E(\bar\Q)$ of points of order dividing 
$p^n$  
and  $T_p(E)$  for the inverse limit of the  $E[p^n]$ with 
respect to
multiplication by  $p$.  For every  $n$,  $E[p^n] \cong 
(\Z/p^n\Z)^2$,  
and so  $T_p(E) \cong \Zp^2$. The action of  $\GQ$ induces 
a representation  
$$
\rhoEp : \GQ  \to  \GL2(\Zp)
$$
such that $\det(\rhoEp) = \e_p$ and
for all but finitely many primes $q$,
\begin{equation}
\label{truees}
\trace(\rhoEp(\Frobq)) = q + 1 - \#(E(\Fq)).
\end{equation}
Composing $\rhoEp$ with the reduction map from $\Zp$ 
to $\Fp$ gives $\rhomod{p}$ of \S\ref{mcec}.

\subsection{Modular representations}
\label{mr}

If  $f$ is an eigenform and $\lambda$ is a prime ideal of 
$\Of$, let
$\Ofl$ denote the completion of $\Of$ at $\lambda$.

\begin{defn}
If $A$ is a ring, a  representation  $\rho : \GQ \to 
\GL2(A)$ is called
{\em modular} if there are an eigenform $f(z) = \fexp$, a 
ring $A'$ 
containing $A$, and a homomorphism  $\iota : \Of \to A'$ 
such that for all but finitely many primes  $q$,
$$
\trace(\rho(\Frobq))  =  \iota(a_q).
$$
\end{defn}

\begin{exs}
(i)  Given an eigenform $f(z) = \fexp$ and a prime ideal 
$\lambda$  of  $\Of$,  Eichler and Shimura (see \S 7.6 of 
\cite{Shimura-red-book}) constructed a representation
$$
\rfl : \GQ  \to  \GL2(\Ofl)
$$
such that $\det(\rfl) = \e_p$ (where $\lambda \cap \Z = 
p\Z$) and for all but 
finitely many primes $q$,
\begin{equation}
\label{fcoeffs}
\trace(\rfl(\Frobq)) = a_q.
\end{equation}
Thus $\rfl$ is modular with $\iota$ taken to be the 
inclusion of 
$\Of$ in $\Ofl$.

(ii)  Suppose $p$ is a prime and $E$ is an elliptic curve 
over $\Q$. If
$E$ is modular, then $\rhoEp$ and $\rhomod{p}$ are modular 
by 
(\ref{truees}), (\ref{es}), and
(\ref{mdlr}).  Conversely, if $\rhoEp$ is modular, then 
it follows from (\ref{truees}) that $E$ is modular.  This 
proves the
following.
\end{exs}
\begin{thm}
\label{mods}
Suppose  $E$  is an elliptic curve over $\Q$.  Then
\begin{center}  
$E$  is modular  $\Leftrightarrow$  $\rhoEp$  is modular 
for every  $p$ 
$\Leftrightarrow$  $\rhoEp$  is modular for one  $p$.
\end{center}
\end{thm}

\begin{rem}
In this language, the Semistable Modular Lifting 
Conjecture says that if $p$
is an odd prime, $E$ is a semistable
elliptic curve over $\Q$, and $\rhomod{p}$ is modular and 
irreducible,
then $\rhoEp$ is modular.
\end{rem}

\subsection{Liftings of Galois representations}
 Fix a prime $p$ and a finite field  $k$ of characteristic 
$p$. Recall that
${\bar k}$ denotes an algebraic closure of $k$.

Given a map $\phi : A \to B$, the induced map from 
$\GL2(A)$ to $\GL2(B)$ will
also be denoted $\phi$.

If $\rho : \GQ \to 
\operatorname{GL}
_2(A)$ is a representation and $A'$ is a ring containing
$A$, we write $\rho \otimes A'$ for the composition of 
$\rho$ with the 
inclusion of $\operatorname{GL}
_2(A)$ in $\operatorname{GL}_2(A')$.

\begin{defn}
If  $\rb : \GQ \to \GL2(k)$ is a representation, we say 
that a representation
$\r : \GQ \to \GL2(A)$  is a {\em lifting} of  $\rb$ (to 
$A$) if $A$ is a 
complete noetherian local $\Zp$-algebra and there exists a 
homomorphism
$\iota : A \to {\bar k}$ such that the diagram
\updiag{A}{\bar k}{\r}{\rb\otimes\bar k}{\iota}
commutes, in the sense that $[\iota \circ \r] = [\rb 
\otimes {\bar k}]$.
\end{defn}

\begin{exs}
(i) If $E$ is an elliptic curve then $\rhoEp$ is a lifting 
of $\rhomod{p}$.

(ii) If $E$ is an elliptic curve, $p$ is a prime, and 
hypotheses (a) and
(b) of Conjecture \ref{ME} hold with an eigenform $f$ and 
prime ideal 
$\lambda$, then $\rfl$ is a lifting of $\rhomod{p}$.
\end{exs}

\subsection{Deformation data}
We will be interested not in all liftings of a given
$\rb$,  but rather in those satisfying
various restrictions.  See Appendix \ref{A1} for the 
definition of the inertia 
groups $I_q \subset \GQ$ associated to primes $q$.  
We say that a representation  $\rho$ of $\GQ$ is {\em 
unramified} at a 
prime  $q$  if  $\rho(I_q) = 1$.  If  $\Sigma$  is a set of 
primes, we say $\r$  is 
{\em unramified outside of  $\Sigma$}  if  $\rho$  is 
unramified at every  $q \notin \Sigma$.

\begin{defn}
By {\em deformation data} we mean a pair  
$$
\D  =  (\Sigma, t)
$$
where $\Sigma$ is a finite set of primes and $t$ is one of 
the words 
{\em ordinary} or {\em flat}.
\end{defn}

If $A$ is a $\Z_p$-algebra, let
$\e_A : \GQ \to \Zp^\times \to A^\times$  be the 
composition of the
cyclotomic character $\e_p$ with the structure map.

\begin{defn}
Given deformation data $\D$, a representation 
$\rho : \GQ \to \GL2(A)$ is {\em type}-$\D$ if $A$ is a 
complete noetherian
local $\Z_p$-algebra, $\det(\rho) = \e_A$, $\rho$  is 
unramified outside of  $\Sigma$,
and $\rho$ is $t$ at $p$ (where $t \in$ \{ordinary, flat\};
see Appendix \ref{A2}).
\end{defn}

\begin{defn}
A representation  $\rb : \GQ \to \operatorname{GL}
_2(k)$  is {\em $\D$-modular} if 
there are an eigenform  $f$ and a prime ideal  $\lambda$ 
of $\Of$
such that $\rfl$ is a type-$\D$ lifting of $\rb$.  
\end{defn}

\begin{rems}
(i) A representation with a type-$\D$ lifting must itself be
type-$\D$. Therefore if a representation is $\D$-modular,
 then it is both  
type-$\D$  and modular.

(ii) Conversely, if $\rb$ is type-$\D$, modular, and 
satisfies (ii) 
of Theorem \ref{BW} below, then $\rb$ is $\D$-modular, by
work of Ribet and 
others (see \cite{Ribetreport}).  This plays an important 
role in Wiles'
 work.
\end{rems}

\subsection{Mazur Conjecture}

\begin{defn}
A representation $\rb : \GQ \to \GL2(k)$  is called {\em 
absolutely 
irreducible} if ${\bar {\rho}} \otimes {\bar k}$ is 
irreducible.
\end{defn}

The following variant of a conjecture of 
Mazur (see Conjecture 18 of \cite{Mazur-Tilouine}; see 
also Conjecture \ref{M}
below) implies the Semistable Modular Lifting Conjecture.

\begin{conj}[Mazur]
\label{M1}
Suppose $p$ is an odd prime, $k$ is a finite field of 
characteristic $p$,
$\D$ is deformation data, and  $\rb : \GQ \to \GL2(k)$  is 
an absolutely irreducible 
$\D$-modular representation.  Then every type-$\D$ lifting 
of $\rb$ to 
the ring of integers of a finite extension of $\Qp$ is 
modular.
\end{conj}

\begin{rem}
Loosely speaking, Conjecture \ref{M1} says that if $\rb$ 
is modular,
then every lifting which ``looks modular'' is modular. 
\end{rem}

\begin{defn}
An elliptic curve $E$ over $\Q$ has {\em good} 
(respectively, {\em bad}\,\RM) 
{\em reduction} at a prime $q$ if $E$ 
is nonsingular (respectively, singular) modulo $q$. An 
elliptic curve $E$ over 
$\Q$ has {\em ordinary} (respectively, {\em 
supersingular}) {\em reduction} 
at $q$ if $E$ has good reduction at $q$ and $E[q]$ has 
(respectively, does 
not have) a subgroup of order $q$ stable under the inertia 
group $I_q$.
\end{defn}

\begin{prop}
\label{cccc}
Conjecture \ref{M1} implies Conjecture \ref{ME}.
\end{prop}

\begin{pf}
Suppose $p$ is an odd prime and $E$ is a semistable 
elliptic curve over $\Q$ 
which satisfies  (a) and (b) of 
Conjecture \ref{ME}. We will apply Conjecture \ref{M1} 
with $\rb = \rhomod{p}$.  
Write $\tau$  for complex conjugation.  Then  $\tau^2 = 
1$,  and by
(\ref{cycl}), $\det(\rhomod{p}(\tau)) = -1$.  Since 
$\rhomod{p}$ is
irreducible and $p$ is odd, a simple linear algebra 
argument now shows 
that $\rhomod{p}$ is absolutely irreducible.

Since $E$ satisfies (b) of Conjecture \ref{ME}, 
$\rhomod{p}$ is modular. Let
\vskip 0.5pc
\begin{itemize}
\item $\Sigma = \{ p \} \cup \{{\text{primes $q$}} : 
{\text{$E$ has bad reduction 
at $q$}} \}$,
\item $t =$ {\em ordinary} if $E$ has ordinary or bad 
reduction at $p$, \\
$t =$ {\em flat} if $E$ has supersingular reduction at $p$,
\item $\D = (\Sigma,t)$.
\end{itemize}
\vskip 3pt
Using the semistability of $E$, one can show that
$\rhoEp$ is a type-$\D$ lifting of $\rhomod{p}$ and
(by combining results of several people; see 
\cite{Ribetreport})
that $\rhomod{p}$ is $\D$-modular.
Conjecture \ref{M1} then says $\rhoEp$ is modular. By 
Theorem \ref{mods},
$E$ is modular.
\end{pf}

\vskip 1pc
\section{Mazur's deformation theory}
\label{universal}
Next we reformulate Conjecture \ref{M1} as a conjecture 
(Conjecture \ref{M})
that the algebras which parametrize liftings and modular 
liftings of a
given representation are isomorphic. It is this form of 
Mazur's conjecture
that Wiles attacks directly.

\subsection{The universal deformation algebra $R$}
 Fix an odd prime $p$, a finite field $k$ of 
characteristic $p$, deformation
data $\D$, and an absolutely irreducible type-$\D$ 
representation
$\rb : \GQ \to \GL2(k)$. Suppose $\O$ is the ring of 
integers of a finite
extension of $\Qp$ with residue field $k$.

\begin{defn}
We say $\r : \GQ \to \GL2(A)$ is a {\em $(\D,\O)$-lifting} 
of $\rb$ if 
$\r$ is type-$\D$, $A$ is a
complete noetherian local $\O$-algebra with residue field 
$k$, and 
the following diagram commutes
\vskip 1pc
\updiag{A}{k}{\r}{\rb}{}
where the vertical map is reduction modulo the maximal 
ideal of $A$.
\end{defn}

\begin{thm}[Mazur-Ramakrishna]
\label{MR}
With $p$, $k$, $\D$, $\rb$, and $\O$ as above,  there are an
$\O$-algebra  $R$  and a  $(\D,\O)$-lifting
$\r_R : \GQ  \to  \GL2(R)$
of  $\rb$, with the property that for every  
$(\D,\O)$-lifting  
$\r$  of  $\rb$ to $A$  there is
a unique $\O$-algebra homomorphism  $\phi_\r : R \to A$ 
such that 
the diagram
\downdiag{R}{A}{\r_R}{\r}{\phi_\r}
commutes.
\end{thm}

This theorem was proved by Mazur \cite{Mazur-def} in the
case when  $\D$  is ordinary and by Ramakrishna 
\cite{Ramakrishna} when  $\D$  
is flat.  Theorem \ref{MR} determines $R$ and $\r_R$ up to 
isomorphism.

\subsection{The universal modular deformation algebra $\T$}
\label{hecke}
 Fix an odd prime $p$, a finite field $k$ of 
characteristic $p$, deformation
data $\D$, and an absolutely irreducible type-$\D$ 
representation
$\rb : \GQ \to \GL2(k)$. Assume $\rb$ is $\D$-modular, and 
fix
an eigenform $f$ and a prime ideal $\l$ of $\Of$ such that
$\rfl$ is a type-$\D$ lifting of $\rb$.  Suppose in 
addition that $\O$ is
the ring of integers of a finite extension of $\Qp$ with 
residue field $k$,
$\Ofl \subseteq \O$, and the diagram
\updiag{\Ofl}{k}{\rfl}{\rb}{}
commutes, where the vertical map is the reduction map.

Under these assumptions $\rfl \otimes \O$ is a 
$(\D,\O)$-lifting of $\rb$,
and Wiles constructs a generalized 
Hecke algebra  $\T$ which has the following properties 
(recall that Hecke
algebras $T(N)$ were defined in \S\ref{modforms}).
\begin{itemize}
\item[$(\T1)$]  \<$\T$  is a complete noetherian local 
$\O$-algebra with residue
field $k$.
\item[$(\T2)$]  There are an integer $N$ divisible only by 
primes in $\Sigma$
and a homomorphism from the Hecke algebra $T(N)$ to $\T$ 
such that $\T$ is
generated over  $\O$  by the images of the Hecke operators  
$T_q$ for primes $q \notin \Sigma$.  By abuse of notation 
we write $T_q$ also 
for its image in $\T$.
\item[$(\T3)$]  There is a $(\D,\O)$-lifting
$$
\r_\T : \GQ \to \GL2(\T)
$$
of $\rb$  with the property that
$\trace(\r_\T(\Frobq))  =  T_q$ for every prime  $q \notin 
\Sigma$.
\item[$(\T4)$] If  $\r$ is modular and is a 
$(\D,\O)$-lifting of $\rb$ to $A$, then
there is a unique $\O$-algebra homomorphism  $\psi_\r : \T 
\to A$  such that 
the diagram
\downdiag{\T}{A}{\r\,{}_\T}{\r}{\psi_\r}
commutes.
\end{itemize}

Since  $\r_\T$   is a $(\D,\O)$-lifting of  $\rb$,  by 
Theorem \ref{MR} there is a
homomorphism
$$  
\j  : R \to \T
$$
such that  $\r_\T$ is isomorphic to $\j \circ \r_R$.  By 
$(\T3)$,  
$\j(\trace(\r_R(\Frobq)))  =  T_q$
for every prime $q \notin \Sigma$, so it follows from 
$(\T2)$ that  $\j$  
is surjective.

\subsection{Mazur Conjecture, revisited}
Conjecture \ref{M1} can be reformulated in the following 
way.

\begin{conj}[Mazur]
\label{M}
Suppose $p$, $k$, $\D$, $\rb$, and $\O$ are as in 
\S\ref{hecke}.
Then the above map $\j : R \to \T$  is an isomorphism.
\end{conj}

Conjecture \ref{M} was stated in \cite{Mazur-Tilouine} 
(Conjecture 18) for
$\D$ ordinary, and Wiles modified the conjecture to 
include the flat case. 

\begin{prop}
Conjecture \ref{M} implies Conjecture \ref{M1}.
\label{imp}
\end{prop}

\begin{pf}
Suppose $\rb : \GQ  \to  \GL2(k)$ is absolutely 
irreducible  and $\D$-modular, 
$A$ is the ring of integers of a finite extension of 
$\Qp$, and
$\r$ is a type-$\D$ lifting of $\rb$ to $A$.
Taking $\O$ to be the ring of integers of a sufficiently 
large finite extension 
of $\Qp$, and extending $\r$ and $\rb$ to $\O$ and its 
residue field, 
respectively, we may assume that $\r$ is a
$(\D,\O)$-lifting of $\rb$.  Assuming Conjecture \ref{M}, 
let
$\psi = \phi_\r \circ \j^{-1} : \T \to A$, with $\phi_\r$ 
as in
Theorem \ref{MR}. By $(\T3)$ and Theorem \ref{MR}, 
$\psi(T_q) = \trace(\r(\Frobq))$ for all
but finitely many $q$.  By \S 3.5 of 
\cite{Shimura-red-book}, given such a 
homomorphism $\psi$ (and viewing $A$ as a subring of $\C$),
there is an eigenform  $\fexp$  where  $a_q = \psi(T_q)$ 
for all but finitely many primes $q$.  Thus $\r$ is modular.
\end{pf}


\section{Wiles' approach to the Mazur Conjecture}
\label{proof}
In this section we sketch the major ideas of Wiles' attack 
on 
Conjecture \ref{M}. The first step (Theorem \ref{WW}), and 
the
key to Wiles' proof, is to reduce Conjecture \ref{M}
to a bound on the order of the cotangent space at a prime 
of $R$.
In \S\ref{selmer} we see that the corresponding tangent 
space is a Selmer
group, and in \S\ref{euler} we outline a general procedure 
due to Kolyvagin  
for bounding sizes of Selmer groups. The input for 
Kolyvagin's method is 
known as an Euler system. The most difficult part of 
Wiles' work 
(\S\ref{geometric}), and the part described as ``not yet 
complete'' in
his December announcement, is his construction of a
suitable Euler system. In \S\ref{fr} we state the results 
announced by Wiles 
(Theorems \ref{BW} and \ref{LW} and Corollary \ref{cm}) 
and explain why 
Theorem \ref{BW} suffices for proving the Semistable 
Taniyama-Shimura 
Conjecture. As an application of Corollary \ref{cm} we 
write down an
infinite family of modular elliptic curves.

 For \S{5} fix $p$, $k$, $\D$, $\rb$, $\O$, $f(z) = 
\fexp$, and $\l$ as
in \S\ref{hecke}.  

By property $(\T4)$ there is a homomorphism  
$$
\pi : \T  \to   \O
$$
such that  $\pi \circ \r_\T$ is isomorphic to  $\rfl 
\otimes \O$. By 
property $(\T2)$ and (\ref{fcoeffs}), $\pi$ satisfies 
$\pi(T_q) = a_q$
for all but finitely many $q$. 

\subsection{Key reduction} 
\label{cot} 
Wiles uses the following generalization of a theorem of 
Mazur, which
says that $\T$ is Gorenstein.

\begin{thm}
\label{gor}
There is a \RM(noncanonical\,\RM)  $\T$-module isomorphism
$$\Hom_\O(\T, \O) \isom  \T.$$
\end{thm}

Let  $\n$  denote the ideal of  $\O$  generated by the 
image 
under the composition
$$
\Hom_\O(\T, \O)   \isom   \T  {\; {\buildrel \pi \over 
\rightarrow} \;} \O
$$
of the element $\pi \in \Hom_\O(\T, \O)$.
The ideal  $\n$  is 
well defined independent of the choice of isomorphism in 
Theorem \ref{gor}.

The map $\pi$ determines distinguished prime ideals of  
$\T$  and  $R$,
$$
\psubt = \ker(\pi), \qquad  \psubr = \ker(\pi \circ \j) = 
\j^{-1}(\psubt).
$$

\begin{thm}[Wiles]
\label{WW}
If  
$$
\#(\psubr/\psubr^2) \leq \#(\O/\n) < \infty,
$$  
then  $\j : R \to \T$  is an isomorphism.
\end{thm}

The proof is entirely 
commutative algebra.  The surjectivity of  $\j$  shows 
that  $\#(\psubr/\psubr^2) \geq \#(\psubt/\psubt^2)$,  
and Wiles proves that  
$\#(\psubt/\psubt^2) \geq \#(\O/\n)$.  Thus if 
$\#(\psubr/\psubr^2) \leq \#(\O/\n)$, then
\begin{equation}
\label{equality}
	\#(\psubr/\psubr^2)  =  \#(\psubt/\psubt^2)  =  \#(\O/\n).
\end{equation}
The first equality in (\ref{equality}) shows that $\j$ 
induces an isomorphism
of tangent spaces. 
Wiles uses the second equality in (\ref{equality}) and 
Theorem \ref{gor}
to deduce that  $\T$  is a local complete 
intersection over  $\O$  (that is, there are $f_1, \ldots, 
f_r \in 
\O[[x_1, \ldots, x_r]]$ such that
$$\T \cong \O[[x_1, \ldots, x_r]]/(f_1, \ldots, f_r)$$  
as $\O$-algebras). 
Wiles then combines these two results to prove 
that  $\j$  is an isomorphism.

\subsection{Selmer groups}
\label{selmer}
In general, if  $M$  is a torsion  $\GQ$-module, a Selmer 
group attached to  $M$  
is a subgroup of the Galois cohomology group  $H^1(\GQ, 
M)$  determined by 
certain ``local conditions'' in the following way.  If  
$q$  is a prime with
decomposition group  $D_q \subset \GQ$, then there is a 
restriction map
$$
\res_q : H^1(\GQ, M)  \to  H^1(D_q, M).
$$
 For a fixed collection of subgroups  $\J = \{J_q 
\subseteq H^1(D_q, M) : 
{\text{$q$ prime}}\}$ depending on the particular problem 
under consideration,
the corresponding Selmer group is
$$
S(M)  =  \bigcap_q \res_q^{-1}(J_q)  \subseteq  H^1(\GQ, M).
$$
Write $H^i(\Q, M)$ for $H^i(\GQ, M)$, and $H^i(\Q_q, M)$ 
for $H^i(D_q, M)$.  

\begin{ex}
The original examples of Selmer groups come from 
elliptic curves.  Fix an elliptic curve  $E$  and a 
positive integer  $m$,  
and take  $M 
= E[m]$,  the subgroup of points in  $E(\bar \Q)$  of 
order dividing  $m$.  
There is a natural inclusion 
\begin{equation}
\label{kummer} 
E(\Q)/mE(\Q)  \hookrightarrow  H^1(\Q, E[m])
\end{equation}
obtained by sending  $x \in E(\Q)$  to the cocycle
$\sigma  \mapsto  \sigma(y) - y$,
where  $y \in E(\bar\Q)$ is any point satisfying  $my = 
x$.  Similarly, for 
every prime  $q$  there is a natural inclusion
$$
E(\Q_q)/mE(\Q_q)  \hookrightarrow  H^1(\Q_q, E[m]).
$$
Define the Selmer group  $S(E[m])$  in this case by taking  
the group  $J_q$  to be the image of  $E(\Q_q)/mE(\Q_q)$  
in  $H^1(\Q_q, E[m])$,
for every  $q$.  
This Selmer group is an important tool in studying the 
arithmetic of  $E$  
because it contains (via (\ref{kummer}))  $E(\Q)/mE(\Q)$.
\end{ex}

\medskip

Retaining the notation from the beginning of \S\ref{proof}, 
let $\m$ denote the maximal ideal of 
$\O$ and fix a positive integer  $n$.
The tangent space  \!$\Hom_\O(\psubr/\psubr^2, 
\O\!/\!\m^n)$  can be identified with 
a Selmer group as follows.   

Let  $\Vn$  be the matrix
algebra $\mathrm{M}_2(\O/\m^n)$, with  $\GQ$  acting via 
the adjoint 
representation $\sigma(B)  =  
\rfl(\sigma)B\rfl(\sigma)^{-1}$.
There is a natural injection 
$$
s : \Hom_\O(\psubr/\psubr^2, \O/\m^n) \hookrightarrow 
H^1(\Q, \Vn)
$$
which is described in Appendix \ref{A3} (see also \S 1.6 
of \cite{Mazur-def}).  Wiles 
defines a collection  $\J = \{J_q \subseteq H^1(\Q_q, 
\Vn)\}$ depending on $\D$. 
Let $S_\D(\Vn)$ denote the associated Selmer group.   
Wiles proves that $s$ induces an isomorphism
$$
\Hom_\O(\psubr/\psubr^2, \O/\m^n)  \isom  S_\D(\Vn).
$$

\subsection{Euler systems}
\label{euler}
We have now reduced the proof of Mazur's conjecture to 
bounding the size of the 
Selmer groups  $S_\D(\Vn)$.  About five years ago 
Kolyvagin \cite{Kolyvagin}, 
building on ideas of his own and of Thaine \cite{Thaine}, 
introduced a 
revolutionary new method for bounding the size of a Selmer 
group.  This new 
machinery, which is crucial for Wiles' proof, is what we 
now describe.
	
Suppose  $M$  is a  $\GQ$-module 
of odd exponent  $m$  and  $\J\! =\! \{J_q\! \subseteq \!
H^1(\Q_q\!,\! M)\!\}$  is a system of subgroups with 
associated Selmer group  
$S(M)$  as in \S\ref{selmer}.  
Let  $\Mhat = \Hom(M, \mn)$, where $\mn$ is the 
group of $m$-th roots of unity.  For every prime $q$, the 
cup product gives 
a nondegenerate Tate pairing
$$
\langle\;,\;\rangle
_q : H^1(\Q_q, M) \times H^1(\Q_q, \Mhat)  \to  H^2(\Q_q, 
\mn) 
\isom \Z/m\Z
$$
(see Chapters VI and VII of \cite{Cassels-Frohlich}).
If  $c \in H^1(\Q, M)$ and $d \in H^1(\Q, \Mhat)$, then
\begin{equation}
\label{reclaw}
\sum_q \langle \res_q(c),\res_q(d)\rangle_q  =  0.
\end{equation}

Suppose that  $\L$  is a finite set of primes.  Let  $\sls 
\subseteq H^1(\Q, 
\Mhat)$ be the Selmer group given by the local conditions  
$\J^* = \{\Jqs 
\subseteq H^1(\Q_q, \Mhat)\}$,  where  
\begin{equation*}
\Jqs = 
\begin{cases}
\text{the orthogonal complement of  $J_q$ under $\langle 
\;,\;\rangle_q$} &
\text{if $q \notin \L$,} \\
H^1(\Q_q, \Mhat) & \text{if $q \in \L$.}
\end{cases}
\end{equation*}
If  $d \in H^1(\Q, \Mhat)$,  define
$$
\kd : \prod_{q \in \L} J_q \to \Z/m\Z
$$
by
$$
\kd((c_q))  =  \sum_{q \in \L} \langle
 c_q,\res_q(d)\rangle _q.
$$
Write  $\res_\L : H^1(\Q, M) \to \prod_{q \in \L} 
H^1(\Q_q, M)$  for the 
product of the restriction maps. By (\ref{reclaw}) and the 
definition of 
$\Jqs$,
if $d \in \sls$, then  $\res_\L(S(M)) \subseteq \ker(\kd)$.
If in addition $\res_\L$  is injective on $S(M)$, then
$$
\#(S(M)) \leq \#\bigl(\bigcap_{d \in \sls} 
\ker(\kd)\bigr).  
$$

The difficulty is to produce enough cohomology classes in 
$\sls$ to show
that the right side of the above inequality is small.  
Following
Kolyvagin, an Euler system is a compatible collection of 
classes 
$\kappa(\L) \in \sls$ for a
large (infinite) collection of sets of primes $\L$. 
Loosely speaking,
compatible means that if $\ell \notin \L$, then 
$\res_\ell(\kappa(\L \cup \{\ell\}))$ is related to 
$\res_\ell(\kappa(\L))$.  
Once an Euler system is given, Kolyvagin has 
an inductive procedure for choosing a set $\L$ such that
\begin{mylist}
\item $res_\L$  is injective on $S(M)$,
\item $\bigcap_{\P \subseteq \L} 
\ker(\theta_{\kappa(\P)})$ can be computed 
in terms of $\kappa(\emptyset)$.
\end{mylist}
(Note that if $\P \subseteq \L$, then 
$S_\P^* \subseteq \sls$, so $\kappa(\P) \in \sls$.)

 For several important Selmer groups 
it is possible to construct Euler systems for which 
Kolyvagin's procedure
produces a set $\L$ actually giving an equality
$$
\#(S(M)) = \#\bigl(\bigcap_{\P \subseteq \L} 
\ker(\theta_{\kappa(\P)})\bigr).
$$
This is what Wiles needs to do for the Selmer group 
$S_\D(\Vn)$.
There are several examples in the literature where this 
kind of argument is 
worked out in some detail.  For the simplest case, where 
the Selmer group in 
question is the ideal class group of a real abelian number 
field 
and the $\kappa(\L)$ are constructed from cyclotomic 
units, see 
\cite{Rubin-Appendix}.  For other cases involving ideal 
class groups and Selmer 
groups of elliptic curves, see \cite{Kolyvagin}, 
\cite{Rubin-Main-Conj}, 
\cite{Rubin-Gauss-Sums}, \cite{Gross}. 

\subsection{Wiles' geometric Euler system}
\label{geometric}  The task now is to construct an
Euler system of cohomology classes with which to 
bound  $\#(S_\D(\Vn))$  using 
Kolyvagin's method.  This is the most technically 
difficult part
of Wiles' proof and is the part of Wiles' work he referred 
to as not yet
complete in his December announcement. We give only 
general remarks about
Wiles' construction.

The first step in the construction is due to Flach 
\cite{Flach}.  He constructed 
classes  $\kappa(\L) \in \sls$  for sets  $\L$  consisting 
of just one prime.  This 
allows one to bound the exponent of  $S_\D(\Vn)$, but not 
its order.

Every Euler system starts with some explicit, concrete 
objects.  Earlier 
examples of Euler systems come from cyclotomic or elliptic 
units, Gauss sums, or 
Heegner points on elliptic curves.  Wiles (following 
Flach) constructs his 
cohomology classes from modular units, i.e., meromorphic 
functions on modular 
curves which are holomorphic and nonzero away from the 
cusps.  More precisely,  
$\kappa(\L)$   comes from an explicit function on the 
modular curve  $X_1(L, N)$,  the 
curve obtained by taking the quotient space of the upper 
half plane by the 
action of the group
$$
\Gamma_1(L, N)  =  \{
\bigl( \begin{smallmatrix} a & b \\ c & d 
\end{smallmatrix} \bigr) 
\in \SL2(\Z) : c \equiv 0  \pmod{LN},  \quad
a \equiv d \equiv 1  \pmod{L}\},
$$
and adjoining the cusps, where  $L = \prod_{\ell \in 
\L}\ell$  and where  $N$  
is the $N$ of $(\T2)$ of \S \ref{hecke}.  The construction 
and study of the
classes $\kappa(\L)$ rely heavily on results of Faltings 
\cite{Faltings1}, 
\cite{Faltings2} and others.

\subsection{Wiles' results}
\label{fr}
Wiles announced two main results (Theorems \ref{BW} and 
\ref{LW} below) in 
the direction of Mazur's conjecture, under two different 
sets of hypotheses 
on the representation $\rb$. Theorem \ref{BW}
implies the Semistable Taniyama-Shimura Conjecture and 
Fermat's Last 
Theorem. Wiles' proof of Theorem \ref{BW} depends on the 
not-yet-complete 
construction of an appropriate Euler system (as in 
\S\ref{geometric}), 
while his proof of Theorem \ref{LW} (though not yet fully 
checked) does 
not. For Theorem \ref{LW}, Wiles bounds the Selmer group 
of \S\ref{selmer} 
without constructing a new Euler system, by
using results from the Iwasawa theory of imaginary 
quadratic fields. (These
results in turn rely on Kolyvagin's method and the Euler 
system of
elliptic units; see \cite{Rubin-Main-Conj}.) 

Since for ease of exposition
we defined modularity of representations in terms of 
$\Gamma_0(N)$ instead of
$\Gamma_1(N)$, the theorems stated below are weaker than 
those announced by
Wiles, but have the same applications to elliptic curves. 
(Note that by
our definition of type-$\D$, if $\rb$ is type-$\D$, then 
$\det(\rb) = \wp$.)

If $\rb$ is a representation of $\GQ$ on a vector space 
$V$, $\Sym2(\rb)$
denotes the representation on the symmetric square of $V$ 
induced by $\rb$.

\begin{thm}[Wiles]
\label{BW}
Suppose $p$, $k$, $\D$, $\rb$, and $\O$ are as in 
\S\ref{hecke} and $\rb$
satisfies the following additional conditions\,\RM:
\begin{mylist}
\item[{\normalshape (i)}]	$\Sym2(\rb)$  is absolutely 
irreducible,
\item[{\normalshape (ii)}] if  $\rb$   
is ramified at  $q$  and  $q \neq p$,  then the 
restriction of  
$\rb$ to $D_q$ is reducible,
\item[{\normalshape (iii)}]	if  $p$ is $3$  or  $5$, 
then for some prime  
$q$,  $p$ divides $\#(\rb(I_q))$.
\end{mylist}
Then $\j : R \to \T$  is an isomorphism.
\end{thm}

Since Theorem \ref{BW} does not yield the full Mazur 
Conjecture (Conjecture \ref{M}) 
for $p = 3$ and $5$, we need to reexamine the arguments of 
\S\ref{overview} to see 
which elliptic curves $E$ can be proved modular using 
Theorem \ref{BW} applied 
to  $\rhomod{3}$  and  $\rhomod{5}$.

Hypothesis (i) of Theorem \ref{BW} will be satisfied if 
the image of  
$\rhomod{p}$  is sufficiently large in  $\GL2(\Fp)$  (for 
example, if  
$\rhomod{p}$  is surjective).  For $p = 3$ and $p = 5$, if 
$\rhomod{p}$
satisfies hypothesis (iii) and is irreducible, then it 
satisfies hypothesis (i).

If $E$ is semistable, $p$ is an odd prime, and 
$\rhomod{p}$ is irreducible and
modular, then $\rhomod{p}$ is $\D$-modular
for some $\D$ (see the proof of Proposition \ref{cccc}) 
and $\rhomod{p}$ satisfies
(ii) and (iii) (use Tate curves; see \S14 of Appendix C of 
\cite{Silverman}).
Therefore by Propositions \ref{imp} and \ref{cccc}, 
Theorem \ref{BW} implies
that the Semistable Modular Lifting Conjecture (Conjecture 
\ref{ME}) holds
for $p = 3$ and for $p = 5$. As shown in \S\ref{overview}, 
the 
Semistable Taniyama-Shimura
Conjecture and Fermat's Last Theorem follow.

\begin{thm}[Wiles]
\label{LW}
Suppose $p$, $k$, $\D$, $\rb$, and $\O$ are as in 
\S\ref{hecke} and $\O$
contains no nontrivial $p$-th roots of unity.  Suppose
also that there are an imaginary quadratic field $F$ of 
discriminant prime to 
$p$ and a character $\chi : \Gal(\bar\Q/F) \to \O^\times$ 
such that the induced
representation $\mathrm{Ind}\chi$ of $\GQ$ is a 
$(\D,\O)$-lifting of $\rb$.
Then $\j : R \to \T$  is an isomorphism.
\end{thm}

\begin{cor}[Wiles]
\label{cm}
Suppose $E$ is an elliptic curve over $\Q$ with complex 
multiplication 
by an imaginary quadratic field $F$ and $p$ is an odd prime
at which $E$ has good reduction.  If  $E'$ is an elliptic 
curve over $\Q$ 
satisfying
\begin{mylist}
\item $E'$ has good reduction at $p$ and
\item $\rhopmod{p}$ is isomorphic to $\rhomod{p}$,
\end{mylist}
then $E'$ is modular.
\end{cor}

\begin{pf*}{Proof of corollary}
Let $\p$ be a prime of $F$ containing $p$, and define
\begin{mylist}
\item $\O$ = the ring of integers of the completion of $F$ 
at $\p$,
\item $k = \O/\p\O$,
\item $\Sigma = \{\text{primes at which $E$ or $E'$ has 
bad reduction}\}
     \cup\{p\}$,
\item $t =$ {\em ordinary} if $E$ has ordinary reduction 
at $p$, \\
$t =$ {\em flat} if $E$ has supersingular reduction at $p$,
\item $\D = (\Sigma, t)$.
\end{mylist}
Let
$$
\chi : \Gal({\bar \Q}/F) \to 
{\mathrm{Aut}}_\O(E[\p^\infty]) \cong \O^\times
$$
be the character giving the action of $\Gal({\bar \Q}/F)$ 
on $E[\p^\infty]$
(where $E[\p^\infty]$ is the group
of points of $E$ killed by the endomorphisms of $E$ which
lie in some power of $\p$).
It is not hard to see that $\rhoEp \otimes \O$ is 
isomorphic to $\mathrm{Ind}\chi$.  

Since $E$ has complex multiplication, it is well known 
that $E$ and
$\rhomod{p}$ are modular.
Since $E$ has good reduction at $p$, it can be shown that 
the discriminant 
of $F$ is prime to $p$ and $\O$ contains no nontrivial 
$p$-th roots of unity.
One can show that all of the hypotheses of 
Theorem \ref{LW} are satisfied with $\rb = \rhomod{p} 
\otimes k$. 
By our assumptions on $E'$, 
$\rhopEp \otimes \O$ is a $(\D,\O)$-lifting of $\rb$, and 
we conclude
(using the same reasoning as in the proofs of Propositions 
\ref{cccc} 
and \ref{imp}) that $\rhopEp$ is modular and hence $E'$ is 
modular.
\end{pf*}

\begin{rems}
(i) The elliptic curves $E'$ of Corollary \ref{cm} are not
semistable.

(ii) Suppose $E$ and $p$ are as in Corollary \ref{cm} and 
$p = 3$ or $5$.  As in Appendix \ref{app2} one can show 
that the elliptic
curves $E'$ over $\Q$ with good reduction at $p$ and with 
$\rhopmod{p}$ isomorphic to $\rhomod{p}$ give infinitely 
many $\C$-isomorphism
classes.
\end{rems}

\begin{ex}
Take $E$ to be the elliptic curve defined by
$$
y^2 = x^3 - x^2 - 3x -1.
$$
Then $E$ has complex multiplication by $\Q(\sqrt{-2})$, 
and $E$ has 
good reduction at $3$.  Define polynomials
\begin{align*}
a_4(t) &= -2430 t^4 - 1512 t^3 - 396 t^2 - 56 t - 3, \\
a_6(t) &= 40824 t^6 + 31104 t^5 + 8370 t^4 + 504 t^3 - 148 
t^2 - 24 t - 1,
\end{align*}
and for each $t \in \Q$ let $E_t$ be the elliptic curve
$$
y^2 = x^3 - x^2 + a_4(t)x +a_6(t)
$$
(note that $E_0 = E$).
It can be shown that for every $t \in \Q$, $\rb_{E_t,3}$ 
is isomorphic to 
$\rhomod{3}$.  If $t \in \Z$ and
$t \equiv 0$ or $1 \pmod{3}$ (or more generally if $t = 
3a/b$
or $t = 3a/b + 1$ with $a$ and $b$ integers
and $b$ not divisible by $3$), then $E_t$ has
good reduction at $3$, for instance because the 
discriminant of $E_t$ is
$$
2^9 (27 t^2 + 10 t + 1)^3 (27 t^2 + 18 t + 1)^3.
$$
Thus for these values of $t$, Corollary \ref{cm} shows 
that $E_t$ is
modular and so is any elliptic curve over $\Q$ isomorphic 
over $\C$ to $E_t$, 
i.e., any elliptic curve over $\Q$ with $j$-invariant 
equal to
$$
\left({{4(27 t^2 + 6 t + 1)(135 t^2 + 54 t + 5)} \over
{(27 t^2 + 10 t + 1)(27 t^2 + 18 t + 1)}}\right)^3.
$$

This explicitly gives infinitely many modular elliptic 
curves
over $\Q$ which are nonisomorphic over $\C$.

(For definitions of complex multiplication, discriminant, 
and $j$-invariant,
see any standard reference on elliptic curves, such as 
\cite{Silverman}.)
\end{ex}


\appendix
\section{Galois groups and Frobenius elements}
\label{A1}
Write  $\GQ = \Gal(\bar\Q/\Q)$.  If  $q$ is a prime number 
and 
$\q$ is a prime ideal dividing  $q$  in the ring of 
integers of  
$\bar\Q$, there is a filtration
$$
\GQ  \supset  D_\q  \supset  I_\q
$$
where the decomposition group $D_\q$ and the inertia group 
$I_\q$ are
defined by
\begin{align*}
D_\q  &=  \{\sigma \in \GQ : \sigma\q = \q\},      \\
I_\q  &=  \{\sigma \in \D_\q : \sigma x \equiv x \pmod{\q} 
  \text{~for 
all algebraic integers~}  x \}.
\end{align*}
There are natural identifications
$$
	D_\q  \cong  \Gal(\bar{\Q}_q/\Q_q),  \qquad   D_\q/I_\q  
\cong  
\Gal(\bar{\F}_q/\F_q),
$$
and  $\Frobmq \in D_\q/I_\q$ denotes the inverse image of 
the 
canonical generator  $x \mapsto x^q$  of  
$\Gal(\bar{\F}_q/\F_q)$.  
If  $\q'$  is another prime ideal above  $q$,  then  $\q' = 
\sigma\q$  for some  $\sigma \in \GQ$  and
$$
	D_{\q'}  =  \sigma D_\q \sigma^{-1}, \qquad  I_{\q'}  =  
\sigma 
I_\q\sigma^{-1},  
\qquad \mathrm{Frob}_{\q'}  =  \sigma \Frobmq \sigma^{-1}.
$$
Since we will care about these objects only up to 
conjugation, 
we will write  $D_q$ and $I_q$.
We will write $\Frobq \in \GQ$ for any representative of a 
$\Frobmq$.  
If  $\r$ is a representation of  $\GQ$ which is 
unramified at  $q$, then $\trace(\r(\Frobq))$ and 
$\det(\r(\Frobq))$ 
are well defined independent of any choices. 

\section{Some details on the proof of Proposition 
\ref{redat3}}
\subsection{\nopunct}
\label{app1}
The modular curve $X_0(15)$ can be viewed as a curve 
defined over $\Q$ in
such a way that the noncusp rational points correspond to 
isomorphism classes (over $\C$)
of pairs $(E',\cc)$ where $E'$ is an elliptic curve over 
$\Q$ and 
$\cc \subset E(\bar\Q)$ is a subgroup of order $15$ stable 
under $\GQ$.
An equation for $X_0(15)$ is  $y^2 = x(x + 3^2)(x - 4^2)$, 
the elliptic curve
discussed in \S\ref{ellipticcurves}. There are eight 
rational points on $X_0(15)$,
four of which are cusps. There are four modular elliptic 
curves, 
corresponding to a modular form for $\Gamma_0(50)$ (see p. 
86 of \cite{antwerp}), 
which lie in the four distinct $\C$-isomorphism classes 
that correspond to the
noncusp rational points on $X_0(15)$.

Therefore every elliptic curve over $\Q$ with a 
$\GQ$-stable subgroup of
order $15$ is modular.  Equivalently, if $E$ is an 
elliptic curve over
$\Q$ and both $\rhomod{3}$ and $\rhomod{5}$ are reducible, 
then $E$ is
modular.

\subsection{\nopunct}
\label{app2}
 Fix a semistable elliptic curve $E$ over $\Q$.  We will 
show that there 
are infinitely many semistable elliptic curves $E'$ over 
$\Q$ such that
\begin{mylist}
\item[(i)]  $\rhopmod{5}$  is isomorphic to  $\rhomod{5}$, 
and
\item[(ii)]  $\rhopmod{3}$  is irreducible.
\end{mylist}
 
Let
$$\Gamma(5) = \{\ABCD \in \SL2(\Z) : \ABCD \equiv 
\bigl(\begin{smallmatrix} 1&0 \\ 0&1 \end{smallmatrix}\bigr)
 \pmod{5}\}.$$
Let $X$ be the twist of the classical modular curve $X(5)$ 
(see
\cite{Shimura-red-book}) by the cocycle induced by 
$\rhomod{5}$, and let
$S$ be the set of cusps of $X$. Then $X$ is a curve 
defined over $\Q$
which has the following properties.
\begin{mylist}
\item The rational points on $X-S$ correspond to 
isomorphism classes of pairs 
$(E',\phi)$ where $E'$ is an elliptic curve over $\Q$ and 
$\phi : E[5] \to E'[5]$ is a $\GQ$-module isomorphism.
\item As a complex manifold $X-S$ is four copies of 
$\H/\Gamma(5)$, 
so each component of $X$ has genus zero.
\end{mylist}
Let $X^0$ be the component of $X$ containing the rational 
point corresponding 
to $(E,\text{identity})$. Then $X^0$ is a curve of genus 
zero defined over $\Q$ 
with a rational point, so it
has infinitely many rational points.  We want to show that 
infinitely many 
of these points correspond to semistable elliptic curves 
$E'$ with 
$\rhopmod{3}$ irreducible.

There is another modular curve $\hat X$ defined over $\Q$, 
with a finite 
set $\hat S$ of cusps, which has the following properties.
\begin{mylist}
\item The rational points on ${\hat X}-{\hat S}$ 
correspond to 
isomorphism classes of triples 
$(E',\phi, \cc)$ where $E'$ is an elliptic curve over $\Q$, 
$\phi : E[5] \to E'[5]$ is a $\GQ$-module isomorphism, and 
$\cc \subset E'[3]$
is a $\GQ$-stable subgroup of order $3$.
\item As a complex manifold ${\hat X}-{\hat S}$ is four 
copies of
$\H/(\Gamma(5) \cap \Gamma_0(3))$.
\item The map that forgets the subgroup $\cc$ induces a  
surjective morphism $\theta : {\hat X} \to X$ defined over 
$\Q$ and of degree 
$[\Gamma(5) : \Gamma(5) \cap \Gamma_0(3)] = 4$.
\end{mylist}

Let ${\hat X}^0$ be the component of ${\hat X}$ which maps 
to $X^0$.
The function field of $X^0$ is $\Q(t)$, and the function 
field of ${\hat X}^0$
is $\Q(t)[x]/f(t,x)$ where $f(t,x) \in \Q(t)[x]$ is 
irreducible and has 
degree $4$ in $x$.  If  $t' \in \Q$ is sufficiently close 
$5$-adically to the
value of $t$ which corresponds to $E$, then the 
corresponding elliptic curve
is semistable at $5$. By the Hilbert Irreducibility 
Theorem we can find a
$t_1 \in \Q$ so that $f(t_1,x)$ is irreducible in $\Q[x]$. 
It is possible
to fix a prime $\ell \neq 5$ such that $f(t_1,x)$ has no 
roots modulo $\ell$. 
If $t' \in \Q$ is sufficiently close $\ell$-adically
to $t_1$, then $f(t',x)$ has no rational roots, and thus 
$t'$ corresponds to a 
rational
point of $X^0$ which is not the image of a rational point 
of ${\hat X}^0$.
Therefore there are infinitely many elliptic curves $E'$ 
over $\Q$ which are
semistable at $5$ and satisfy
\begin{mylist}
\item[(i)]  $E'[5] \cong E[5]$ as $\GQ$-modules, and
\item[(ii)]  $E'[3]$ has no subgroup of order $3$ stable 
under $\GQ$.
\end{mylist}
It follows from (i) and the semistability of $E$ that $E'$ 
is semistable
at all primes $q \neq 5$, and thus $E'$ is semistable. We 
therefore have
infinitely many semistable elliptic curves $E'$ which 
satisfy the desired 
conditions.

\section{Representation types}
\label{A2}
Suppose $A$ is a complete noetherian local $\Z_p$-algebra 
and 
$\rho : \GQ \to \GL2(A)$ is a representation. Write 
$\rho\mid_{D_p}$ for 
the restriction of $\rho$ to the decomposition group $D_p$.
We say  $\rho$ is

\begin{mylist}
\item {\em ordinary} at $p$ if
$\rho\mid_{D_p}$ is (after a change of basis, if necessary) 
of the form  $\bigl( \begin{smallmatrix} * & * \\ 0 & \chi 
\end{smallmatrix} \bigr)$
where $\chi$ is unramified and the  * are functions from 
$D_p$ to $A$;
\item {\em flat} at $p$ if $\rho$ is not ordinary, and
for every ideal  $\a$  of finite index in  $A$,  the 
reduction of
$\rho\mid_{D_p}$ modulo $\a$  is 
the representation associated to the ${\bar \Q}_p$-points 
of a finite flat group 
scheme over $\Zp$.
\end{mylist}

\section{Selmer groups}
\label{A3}
With notation as in \S\ref{proof} (see especially 
\S\ref{selmer}), define  
$$
\O_n = \O[\eps]/(\eps^2, \m^n)
$$
where $\eps$ is an indeterminate.  Then  $v  \mapsto  1 + 
\eps v$ defines 
an isomorphism
\begin{equation}
\label{map}
\Vn  \isom  \{\delta \in \GL2(\O_n) : \delta \equiv 1 
\pmod{\eps}\}.
\end{equation}

 For every $\alpha \in \Hom_\O(\psubr/\psubr^2, \O/\m^n)$ 
there is a 
unique $\O$-algebra homomorphism  $\y_\alpha : R \to \O_n$ 
 whose
restriction to $\psubr$ is $\eps \alpha$.
Composing with the representation $\r_R$ of Theorem 
\ref{MR} gives a
$(\D,\O)$-lifting $\r_\alpha = \y_\alpha \circ \r_R$  of 
$\rb$ to $\O_n$.
(In particular $\r_0$ denotes the $(\D,\O)$-lifting 
obtained when $\alpha = 0$.)
Define a one-cocycle  $c_\alpha$  on  $\GQ$  by
$$
c_\alpha(g)  =  \r_\alpha(g)\r_0(g)^{-1}.
$$
Since  $\r_\alpha \equiv \r_0 \pmod{\eps}$, using 
(\ref{map}) we can view  
$c_\alpha \in H^1(\Q, \Vn)$.  This defines a homomorphism
$$
s : \Hom_\O(\psubr/\psubr^2, \O/\m^n) \to H^1(\Q, \Vn),
$$
and it is not difficult to see that $s$ is injective.  The 
fact 
that  $\r_0$ and $\r_\alpha$  are type-$\D$ gives 
information 
about the restrictions  $\res_q(c_\alpha)$  for various 
primes  $q$,
and using this information Wiles defines a Selmer group  
$S_\D(\Vn) 
\subset H^1(\Q, \Vn)$ and verifies that $s$ is an 
isomorphism 
onto  $S_\D(\Vn)$.


\end{document}